\documentclass[a4paper]{article}
\usepackage[utf8]{inputenc}
\usepackage{amsfonts}
\usepackage{amsmath}
\usepackage{amssymb}
\usepackage{amsthm}
\usepackage{enumerate}
\usepackage[mathcal]{euscript}
\usepackage{graphicx}
\newtheorem{thm}{Theorem}[section]
\newtheorem{cor}{Corollary}
\newtheorem{prop}[thm]{Proposition}

\theoremstyle{remark}
\newtheorem*{rem}{Remark}

\theoremstyle{definition}
\newtheorem{defin}{Definition}[section]

\newtheorem{exmp}{Example}[section]
\numberwithin{equation}{section}

\title{On the classification of singular flat structures on surfaces}
\author{Ousama Malouf}
\date{}
\begin{document}

\maketitle

\begin{abstract}
We study in this work flat surfaces with conical singularities, that is, surfaces provided with a flat structure with conical singular points. Finding good parameters for these surfaces in the general case is an open question. We give an answer to this question in the case of flat structures on pairs of pants with one singular point. The question of decomposability of an arbitrary flat surface into flat pairs of pants is discussed.
\end{abstract}


\section{Introduction\label{intro}}
In this article, we consider flat surfaces with a finite number of conical singularities, that is, surfaces provided with a flat structure with conical singular points. A conical singular point has a total angle different from $2\pi$.

Finding good parameters for these surfaces is still an open question (one should mention here that there are well-known good parameters for the subset of flat surfaces defined by quadratic differentials on Riemann surfaces). We treat here the question of classifying flat pairs of pants with one singularity. The decomposition of surfaces into pairs of pants is a common practice to study various structures on surfaces (see for instance \cite{pre05560297}). The idea, which is usually attributed to Grothendieck, and which is developed by Feng Luo in \cite{pre05560297}, is to provide building blocks (generators) which permit to reconstruct any surface endowed with the studied structure following some specific rules (relations). It should be noted however that there exist flat surfaces which are not decomposable into flat pairs of pants by disjoint simple closed  geodesics (for an example of flat surface of genus $3$ with one singularity, see Example \ref{nondecomp}). But, it is possible to decompose such a surface into pairs of pants if we sacrifice some rules of decomposition, such as the simplicity of geodesics in the decomposition.

In Section \ref{prelim}, we give a definition of flat surface, we recall the formula of Gauss-Bonnet and we give some examples. In Section \ref{class}, we classify flat pairs of pants with one singularity and we study their space of parameters. In Section \ref{teich}, we present the Teichm\"uller space of flat pairs of pants with one singularity and we study its topology. In Section \ref{exmps}, we discuss the decomposability of a surface into pairs of pants.

This work uses ideas from \cite{1149.57029,1020.57003} in which the authors consider the similar question for hyperbolic surfaces with singular points.

\section{Preliminaries \label{prelim}}
\begin{defin}
Let $M_{g, k}$ (denoted by $M$ if there is no confusion) be a surface of genus $g$ with $k$ boundary components, provided with a flat metric $d (., .)$ with finitely many conical singularities $\varSigma= \{s_1,\ldots, s_n\}$. We assume all the boundary components of $M$ are closed curves (homeomorphic to circles). We denote by $\partial M$ the boundary of $M$ and by $\text{Int}(M)=M\backslash\partial M$ the interior of $M$. More precisely:
\begin{itemize}
 \item For each point $x\in \text{Int}(M) \backslash\varSigma$ there is a neighborhood $U(x)$  isometric to a disc in the Euclidean plane.

\item  For each point $s\in\varSigma \cap \text{Int}(M)$ there is a neighborhood isometric to a Euclidean cone of total angle $0<\theta<+\infty$, with $\theta\neq2\pi$.

\item  For each point $x\in\partial M\backslash\varSigma$ there is a neighborhood isometric to a Euclidean sector with angle measure $\pi$ at the image of $x$.

\item  For each point $s\in \varSigma\cap\partial(M)$ there is a  neighborhood isometric to a Euclidean sector of angle $0<\theta<+\infty$, with $\theta\neq\pi$.

\end{itemize}

The metric $d(.,.)$ will be called a \emph{flat structure} on $M$, and $M$ will be called a \emph{flat surface} (\emph{with boundary} if it exists).

\end{defin}
 
\begin{defin}
 The curvature $\kappa$ at a conical singularity $s$ of total angle $\theta$ is $\kappa=2\pi-\theta$. The curvature at a singular point on the boundary of angle $\theta$ is $\kappa=\pi-\theta$. 
\end{defin}
This definition of curvature is motivated by the following formula:

\begin{prop} [Gauss-Bonnet formula for closed surfaces {\cite[p. 113]{0669.53001}}, {\cite[p. 190]{0146.44103}} or {\cite[p. 85-86]{0611.53035}} for a proof]
Let $M$ be a closed flat surface of genus $g$ with $n$ conical singularities. Let $\theta_i,\,i=1,\ldots,n$ be the total angles at the singularities. The formula of Gauss-Bonnet which connects the number of singularities with their total angles and the genus is:

\begin{equation}
\sum_{i=1}^{n} (2\pi - \theta_i) = (4-4g)\pi 
\label{gb}
\end{equation} 

\end{prop}
The Euler characteristic of this surface is given by $\chi(M)=2-2g$, and the formula of Gauss-Bonnet in terms of the Euler characteristic is:
\[\sum_{i=1}^{n} (2\pi - \theta_i) = 2\pi\chi(M) \,.\]
\begin{cor}[Gauss-Bonnet formula for a disc with $b$ holes]
Let $M$ be a flat surface of genus $0$ with boundary, with $n$ singularities in the interior and $m$ singularities on the boundary, and let $b$ be the number of boundary components. Let $\theta_i,\,i=1,\ldots,n$ be the total angles of conical singularities in the interior, and let $\tau_j,\,j=1,\ldots,m$ be the total angles of singularities on the boundary. Then,
\begin{equation}
 \sum_{i=1}^{n} (2\pi - \theta_i) + \sum_{j=1}^{m} (\pi - \tau_j)= (4-2b)\pi \, .
\label{gbc}
\end{equation}
\end{cor}
\begin{proof}
 By taking the double of $M$ (in this operation, each singularity on the boundary is glued to its copy) we obtain a closed flat surface of genus $g=b-1$ with $2n+m$ singularities. By the formula of Gauss-Bonnet (\ref{gb}) we get the result.
\end{proof}
More generally, the Euler characteristic of a surface $M$ of genus $g$ with $b$ boundary components is given by $\chi(M)=2-2g-b$, and the formula of Gauss-Bonnet in terms of the Euler characteristic is:
\[\sum_{i=1}^{n} (2\pi - \theta_i) + \sum_{j=1}^{m} (\pi - \tau_j)= 2\pi\chi(M) \,.\]

\begin{exmp}
Let $M$ be a pair of pants (a surface homeomorphic to a disc with two holes), equipped with a flat structure with one conical singularity $s$ in its interior and such that each component of the boundary $\partial M=c_1\cup c_2 \cup c_3$ is geodesic without singularities.

From $s$ we take three geodesic segments $d_i, i=1, 2,3$ which realize the distances between $s$ and $c_i$ respectively, Figure \ref{pantc}. We denote by $l_i, r_i$ the lengths of $c_i, d_i$ respectively. We cut and open $M$ along the geodesic segments $d_i$. Then we obtain a connected surface $P$ which can be decomposed into three rectangles $R_i,\,i=1,2,3$ whose lengths of sides are $l_i, r_i$ respectively, and one triangle whose sides have lengths $l_i$,  as shown in Figure \ref{pantc1}.

Following the inverse procedure, we consider a triangle $T$ of side lengths $l_i$, and three rectangles $R_i$ each of which shares a side with $T$ and the other side being of length $r_i$ as shown in Figure \ref{pantc1}.

Now, by identifying the sides of equal length $r_i$,  we obtain a flat pair of pants with one conical singularity. Note that the total angle of the singularity is $\theta=4\pi$.

This example will be studied later in Section \ref{class}.
\end{exmp}
\begin{exmp}
We cannot build a flat pair of pants with no singularities at all. Indeed, if we assume that such a pair of pants exists, by the formula (\ref{gbc}) we find $0=-2\pi$ which is impossible. 
\end{exmp}

Finally, by gluing along their boundaries $2g-2$ flat pair of pants, each one with a single singularity, we obtain a closed flat surface of genus $g$ with $2g-2$ conical singularities.

\section{The geometry of a flat pair of pants with a conical singularity \label{class}}

\begin{figure}
 \centering
 \includegraphics[width=0.35\textwidth,bb=0 0 457 457]{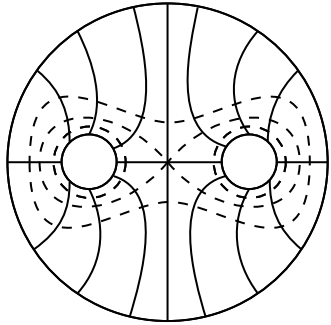}
\caption{The two transverse measured foliations on the pair of pants induce a flat structure with one singular point\label{teichf}} 
\end{figure}

In this section we assume that $M$ is a flat (singular) surface with boundary homeomorphic to a disc with two holes. We assume further that $M$ has only one conical singularity $s$. We call such a surface a \emph{flat pair of pants with one singularity} (or just \emph{flat pair of pants}). An example of such structure is obtained by taking the metric associated to two transverse measured foliations with one singularity on $M$ (see Figure \ref{teichf}). It is clear from this picture that the parameter space for the transverse measures of these two transverse foliations has dimension four, and it produces a space of flat structures of dimension four on the pair of pants, with one singular point. It will follow from the discussion below (see the remark after Definition \ref{teichm}) that there are examples of flat structures on $M$ which do not arise from pairs of transverse measured foliations. The aim of this section is to find a set of real parameters which determine the geometry of $M$.

\begin{figure}
 \centering
 \includegraphics[width=0.25\textwidth,bb=0 0 430 495]{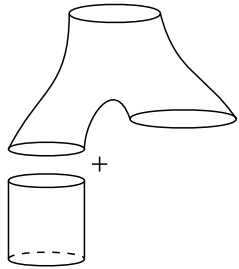}
 \caption{We can change the isometry type of a flat pair of pants without changing the boundary curve lengths\label{cylind}}
\end{figure}

The first natural guess is that the lengths of boundary geodesics might be good parameters, in analogy with the case of hyperbolic pairs of pants (described for instance in \cite{0731.57001}). This is easily seen to be false, as one can glue Euclidean cylinders to boundaries of flat pairs of pants, changing the isometry type, without changing the boundary component length (see Figure \ref{cylind}). Thus, the parameter space for flat structures on pair of pants is more complicated than the parameter space for hyperbolic structure on such surfaces.

Let  $\partial M = c_1 \cup c_2 \cup c_3$ denote the boundary and denote by $d(.,.)$ the distance function on $M$. Let $d_i$ be a geodesic segment which realizes the distance $d (s, c_i)$ between the singularity $s$ and the boundary component $c_i$. This segment intersects the boundary component with a right angle. We denote by $l_i$ the length of $c_i$ and by $r_i$ the length of $d_i$. Obviously, the parameters $r_i,\,i=1,2,3$ depend on the position of $s$. We have the following proposition.

\begin{prop}
  The real parameters $l_i, r_i, i=1, 2,3$ determine a unique flat pair of pants $M$ with one singularity up to isometry, where the $l_i$ are the lengths of the boundary components $c_i$, and the $r_i$ are the lengths of the geodesic segments between $s$ and the $c_i$.
\label{params}
\end{prop}
\begin{proof}
Every boundary component $c_i,\, i\in\{1,2,3\}$ is geodesic, even if the singularity is on a boundary component, since by the formula of Gauss-Bonnet, the angle at $s$ will be $3\pi$. By our assumption on the uniqueness of the singularity, $d_i$ does not meet another singularity on the component $c_i$ of the boundary, and so, since $c_i$ is geodesic, $d_i$ is orthogonal to $c_i$. For the same reasons, $d_i$ does not share any of its interior points with boundary components.

\textit{First case:} We assume that $s\in Int(M)$. This implies that the lengths of $d_i$ are different from zero, $r_i\neq 0, i=1,2,3$. We cut $M$ along $d_1\cup d_2 \cup d_3$, see Figure \ref{pantc}, to obtain a polygonal flat  surface $P$ homeomorphic to a disc without singularities in the interior, see Figure \ref{pantc1}.
\begin{figure}
 \centering
 \includegraphics[width=0.35\textwidth, bb=0 0 430 259]{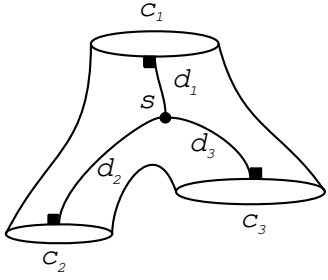}
\caption{\label{pantc}Flat pair of pants with one singularity in the interior}
\end{figure}
\begin{figure}
 \centering
 \includegraphics[width=0.4\textwidth,bb=0 0 517 385]{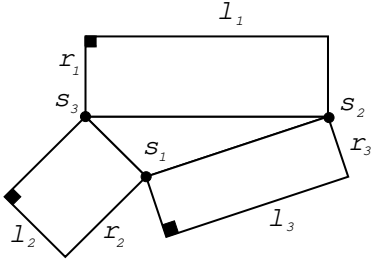}
\caption{\label{pantc1}After cut. $l_i,\,r_i$ are the lengths of $c_i,\,d_i$ respectively}
\end{figure}

None of the parameters $l_i, i=1, 2, 3$ can be zero, because if one of them is zero, the surface $M$ looses one of its boundary components. This fact implies that we always have three copies of $s$ after the cut. Let us denote by $s_i, i=1,2,3$ the copy which lies between $d_{(i+1)\mod 3}$ and $d_{(i+2)\mod 3}$, Figure \ref{pantc1}. We can draw in $P$ the geodesics between these copies. These geodesic segments bound a triangle $T=(s_1 s_2 s_3)$ , Figure \ref{pantat}, which can be degenerate, with none of its angles greater than $\pi$. For instance, in Figure \ref{pantap} the angle $\measuredangle(s_2 s_1 s_3)= \pi$. This angle cannot be greater than $\pi$, because this implies that $d_1$, along which the cut was made, was not a geodesic segment which represents the distance between $s$ and $c_1$. This contradicts the assumption. 
\begin{figure}
 \centering
 \includegraphics[width=0.35\textwidth,bb=0 0 430 401]{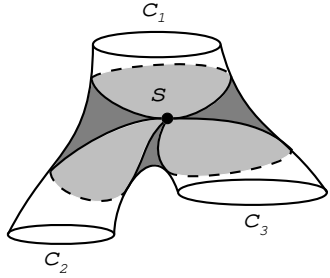}
 \caption{\label{pantat}Here we see the triangle on the pair of pants}
\end{figure}

\begin{figure}
 \centering
 \includegraphics[width=0.4\textwidth,bb=0 0 448 380]{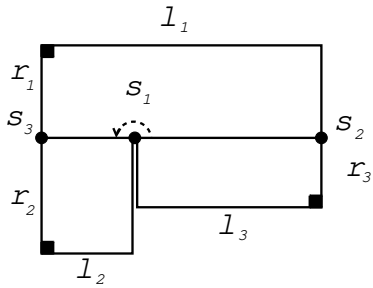}
\caption{\label{pantap}The triangle is degenerate}
\end{figure}

We recognize easily, Figure \ref{pantc1}, in addition to $T$, three rectangles $R_i,\, i=1,2,3$. A rectangle $R_i,\,i\in\{1,2,3\}$ is bounded by a side of $T$, $c_i$ and the two copies of $d_i$. $P$ is composed of $T$ and $R_i,\, i=1,2,3$.

A rectangle is uniquely determined up to isometry by its length and height, and a triangle is uniquely determined up to isometry by its lengths of sides. Then, given the parameters $l_i,r_i,\,i=1,2,3$, where the $l_i,\,i=1,2,3$ satisfy the triangle inequalities, we can construct a unique pair of pants with one singularity by gluing together the three rectangles $R_i$ and the triangle $T$ given by the parameters $l_i, r_i,\, i=1, 2,3$.

The triangle could be degenerate, and this case lies on the boundary of the space of triangles. This happens only when
\begin{equation}
 l_i=l_{(i+1)\mod 3} + l_{(i+2)\mod 3} 
\label{degen}
\end{equation} 
for some $i \in \{1,2,3\}$. Figure \ref{pantap}.

\textit{Second case:} We assume that $s\in \partial M$. Then at least one of $r_i,\, i=1,2,3$  equals zero.
More precisely, $s$ belongs to only one $c_i,\,i\in\{1,2,3\}$, and so only one $r_i$ equals zero, because if $s$ belongs to more than one boundary component, it does not have a neighborhood homeomorphic to a disc. In such a case we call this a \emph{degenerate pair of pants}. We will outline these cases later.
\begin{figure}
 \centering
 \includegraphics[width=0.35\textwidth,bb=0 0 430 401]{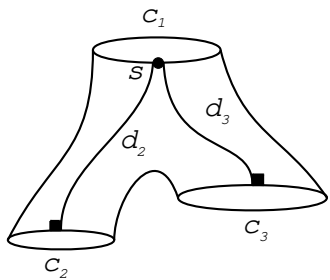}
 \caption{\label{pantabord}Singularity on the boundary} 
\end{figure}
 
Without loss of generality, assume $s\in c_1$, as in Figure \ref{pantabord}. Cut along $d_2 \cup d_3$ ($d_1$ being reduced to $s$) to obtain a polygonal flat surface $P$ homeomorphic to a disc without singularities in the interior. For the same reason as before, we have three copies of $s$. The geodesic segment between $s_2,s_3$ is $c_1$ itself since $c_1$ is geodesic.

The geodesic segments between the copies of $s$ bound a triangle $T=(s_1 s_2 s_3)$ which can be degenerate, with none of its angles greater than $\pi$ (for the same reason as in the first case). It is easy to recognize, in addition to $T$, two rectangles $R_i,\, i=2,3$. Figure \ref{hepta}.
\begin{figure}
 \centering
 \includegraphics[width=0.4\textwidth,bb=0 0 517 274]{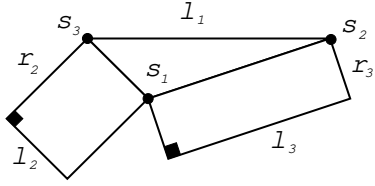}
\caption{\label{hepta}The resulting heptagon after cut} 
\end{figure}

Knowing the parameters $l_i,r_i\,,i=1,2,3$ we can determine if the singularity is on the boundary and on which component. The rectangles $R_i$ and the triangle $T$ are uniquely determined by the parameters and so is the flat structure.

If Relation (\ref{degen}) is satisfied, the triangle is degenerate. In this case we notice that only one $r_i,\, i\in\{1,2,3\}$ that correspond to the boundaries on the right hand side of this identity can take the value zero. The other cases are degenerate.

We proved that to each flat structure on $M$ with one singular point we may associate unique set of values $l_i, r_i\,, i=1,2,3$. Conversely, to each choice of values, we may associate a unique flat structure on $M$.

We deduce that each flat pair of pants $M$ is uniquely determined, up to isometry, by the parameters $l_i$ end $r_i$.  And vice versa.

\end{proof}
Let us give the limits of parameters in view of the last proof. The length parameters can take  any values satisfying these conditions:
\begin{itemize}
 \item $0<l_i<\infty,\,i=1,2,3$,
\item $l_i\leq l_{(i+1)\mod 3}+l_{(i+2)\mod 3} \, ,\, i=1,2,3$,
\item $0\leq r_i < \infty,\,i=1,2,3$,
\item $0<r_i+r_{(i+1)\mod 3},\,i=1,2,3$
\item If $l_i=l_{(i+1)\mod 3}+l_{(i+2)\mod 3} \, , \text{for some }\, i\in\{1,2,3\}$ then $0<r_i<\infty$ 
\end{itemize}

\begin{rem}
Let $M$ be a flat surface with one singularity $s$, which is homeomorphic to a disc with $n$ holes, $n\geq 3$. We denote by $c_i$ the boundary components of $M$ and by $d_i$ the geodesic segments from $s$ to $c_i$. If we set $l_i$ the lengths of $c_i$ and $r_i$ the lengths of $d_i$, then the parameters $l_i, r_i$ do not characterize $M$ in the sense of Proposition \ref{params}. This follows because if we cut $M$ along $d_1 \cup d_2 \cup \dots \cup d_n$ then instead of a triangle $(s_1 s_2 s_3)$ we get a $(n+1)$-gon which, of course, is not uniquely determined  by the lengths of its edges.
 
\end{rem}

We return to the case where $M$ is a pair of pants and $\partial M = c_1 \cup c_2 \cup c_3$. By the previous discussion, where we saw how a flat pair of pants with one singularity can be cut into pieces, three rectangles and one triangle, we find it convenient to introduce these terms: We say that \emph{the triangle is degenerate} if the condition:
\[l_i=l_{(i+1)\mod 3}+l_{(i+2)\mod 3} \, , \text{for some }\, i\in\{1,2,3\}\]
is satisfied. If not, we say that \emph{the triangle is non-degenerate}. We say that \emph{the rectangle} $R_i\,, i\in\{1,2,3\}$ \emph{is degenerate} if $r_i=0$ for some $i\in\{1,2,3\}$. If not, we say that  \emph{the rectangle} $R_i$ is \emph{non-degenerate}.

\begin{rem}[Degenerate cases]
 A degenerate pair of pants appears when the singularity does not have a neighborhood homeomorphic to a disc in the Euclidean plane. From the previous proof, we can conclude the following degenerate cases:
\begin{itemize}
 \item In the case when the triangle is non-degenerate, if two or three rectangles $R_i,\, i\in\{1,2,3\}$ are degenerate we have a degenerate pair of pants. A neighborhood of the singularity is similar to one of the forms in Figure \ref{neibh}.

\begin{figure}
\centering
\begin{tabular}{c c c c c}
 \includegraphics[width=0.25\textwidth,bb=0 0 343 258]{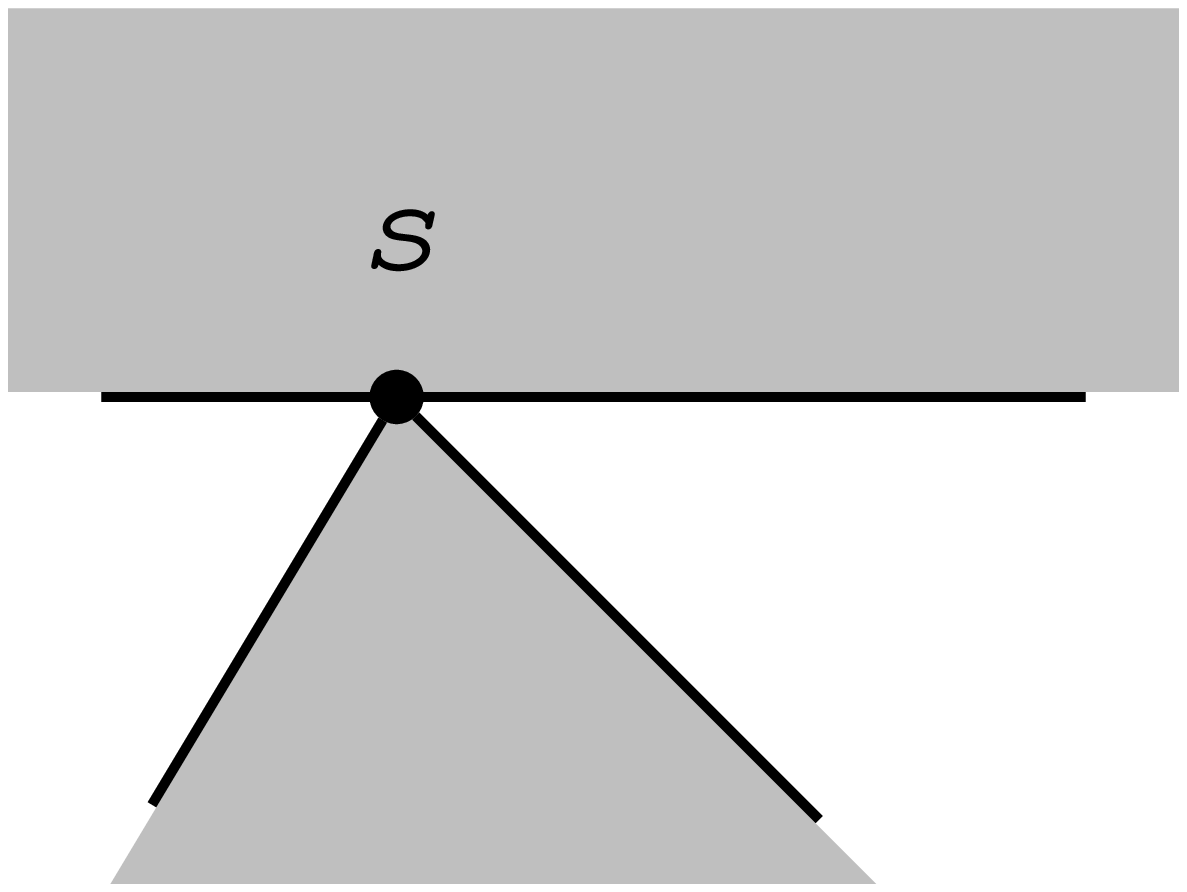}
 & & & &
 \includegraphics[width=0.20\textwidth,bb=0 0 318 316]{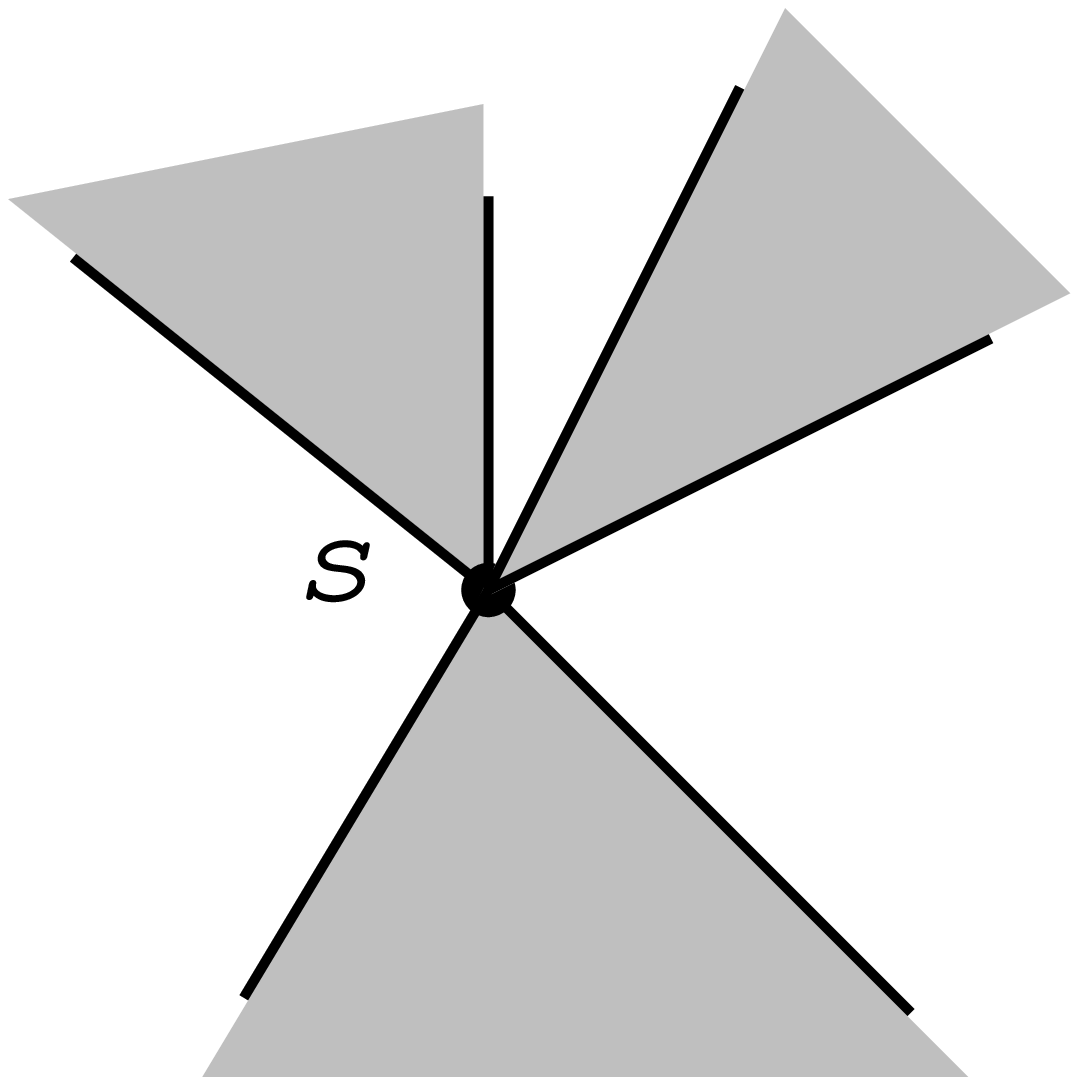}
\\ \\
Two rectangles are degenerate
& & & &
Three rectangles are degenerate
\\
 \end{tabular} 
\caption{\label{neibh}Neighborhood of singularity}
\end{figure} 

\item In the case when the triangle is degenerate, that is, for some $i\in\{1,2,3\}$ identity (\ref{degen}) is satisfied, if both rectangles $R_{(i+1)\mod 3},R_{(i+2)\mod 3}$ are degenerate we have a degenerate pair of pants. A neighborhood of the singularity is similar to Figure \ref{neibh2}.
\begin{figure}
 \centering
 \includegraphics[width=0.20\textwidth,bb=0 0 258 173]{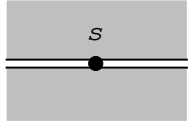}
\caption{\label{neibh2}The triangle and two rectangles are degenerate}
\end{figure}
\end{itemize}

We should finally pay attention that in this case when $R_i$ is degenerate, the resulting surface is not a pair of pants since it has four holes. Figure \ref{imposib}.
\begin{figure}
 \centering
 \includegraphics[width=0.25\textwidth,bb=0 0 288 227]{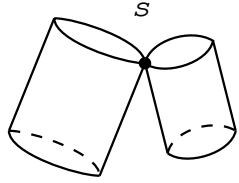}
\caption{\label{imposib}Impossible case}
\end{figure}
\end{rem}

\begin{figure}
 \centering
 \includegraphics[width=0.35\textwidth,bb=0 0 430 401]{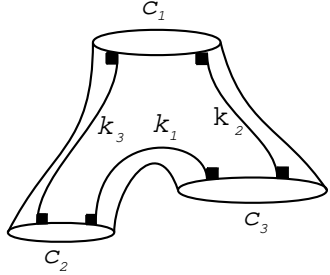}
\caption{\label{pantapa}The geodesic $k_i$ represents the distance $a_i$ between $c_{(i+1)\mod 3}$ and $c_{(i+2)\mod 3}$}
\end{figure}

After the previous discussion, we will introduce new parameters for flat pairs of pants. These parameters seem to be more convenient. For this, let now $k_i$ be a geodesic segment which realizes the distance between $c_{(i+1)\mod 3}$ and $c_{(i+2)\mod 3},\, i=1, 2, 3$, see Figure \ref{pantapa}.  Denote by $l_i$ the length of $c_i$ and by $a_i$ the length of $k_i$.

\begin{defin}
The six non-negative parameters $l_i, a_i, i=1,2,3$ will be called the \emph{distance parameters} of $M$. 
\end{defin}

We have the following theorem.

\begin{prop}
\label{params2}
The distance parameters $l_i, a_i, i=1,2,3$ determine a unique flat pair of pants $M$ with one singularity. 
\end{prop}
\begin{proof}
By the proof of Proposition \ref{params} we see that flat pairs of pants with one singularity are uniquely determined by the parameters $l_i, r_i, i=1,2,3$, and correspond to one of the following cases:
\begin{enumerate}
 \item If the triangle is non-degenerate, we distinguish two cases. 

\begin{enumerate}
 \item If all the rectangles $R_i$ are non-degenerate, it is easy, after the cut, to see that geodesic segments between boundaries pass all by the singularity. Figure \ref{nondegen}. 
\begin{figure}
 \centering
 \includegraphics[width=0.4\textwidth,bb=0 0 536 382]{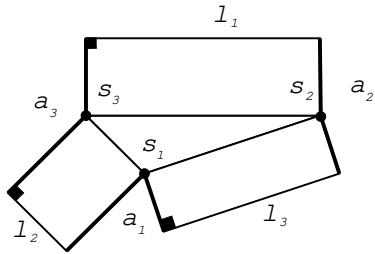}
\caption{\label{nondegen}After cut, $l_i,\,a_i$ are the lengths of $c_i,\,k_i$ respectively}
\end{figure}

So we have:
\begin{equation}
r_i+r_{(i+1) \mod 3}=a_{(i+2)\mod 3}, \, i=1,2,3 \,.
\label{rig} 
\end{equation}

These relations prove that knowing $l_i, a_i,\, i=1,2,3$ we can determine $ r_i,\, i=1, 2, 3 $, and so, all the parameters needed to determine the flat pair of pants in this case, by Proposition \ref{params}.

\item If only one rectangle is degenerate then the pair of pants can be cut into two rectangles and one non-degenerate triangle as in Figure \ref{hepta}. By the same relations (\ref{rig}) we can determine the flat pair of pants. The singularity, in this case, is on the boundary component which corresponds to the degenerate rectangle.
\end{enumerate}

\item If the triangle is degenerate, we also distinguish two cases.
\begin{enumerate}
\item If all the rectangles $R_i,\,i=1,2,3$ are non-degenerate, we see that the relations (\ref{rig}) also hold here and so the flat pair of pants is well determined, in this case, by the distance parameters. Figure \ref{pantap}.

\item If one of $R_{(i+1)\mod 3}, R_{(i+2)\mod 3}$ is degenerate then the pair of pants can be cut into two rectangles. Again, it is determined by relations (\ref{rig}). Figure \ref{2rect}.
\begin{figure}
 \centering
 \includegraphics[width=0.4\textwidth,bb=0 0 456 313]{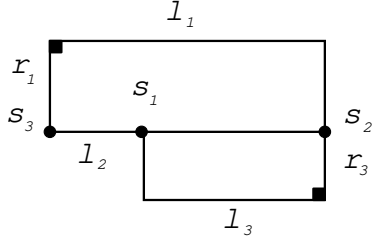}
 \caption{\label{2rect}Permitted gluing of two rectangles}
\end{figure}
\end{enumerate}
\end{enumerate}

\end{proof}
Here we give the limits of the new parameters in view of the last proof. The distance parameters can take  any values in the limit of these conditions:
\begin{enumerate}[1) ]
\item $0<l_i<\infty,\,i=1,2,3$,
\item $l_i\leq l_{(i+1)\mod 3} + l_{(i+2)\mod 3} \, ,\, i=1,2,3$,
\item $0< a_i < \infty,\,i=1,2,3$,
\item $a_i\leq a_{(i+1)\mod 3} + a_{(i+2)\mod 3} \, ,\, i=1,2,3$, (Using Relations (\ref{rig})),
\item If $l_i=l_{(i+1)\mod 3}  + l_{(i+2)\mod 3}$ , for some $i\in\{1,2,3\}$, then \\ $a_i<a_{(i+1)\mod3}+a_{(i+2)\mod3}$ . (Using Relations (\ref{rig})). 
\end{enumerate}

It is easy to see that these conditions are equivalent to preceding ones. In addition, it is interesting to see that these conditions are equivalent to the fact that we have two Euclidean triangles (first four conditions) for which a degenerate case of one triangle prevents a degenerate case of the other (fifth condition).

We denote by $\mathcal{C}$ the set of all flat structures defined by the parameters $l_i,\, a_i,\, i=1,2,3$ under the previous conditions, and by $\mathcal{B}$ the set of all 6-tuples $(l_1,l_2,l_3,a_1,a_2,a_3)\in \mathbb{R}^6$ which satisfy these conditions. We can define a one-to-one mapping $\Phi: \mathcal{C}\rightarrow\mathcal{B}$ such that to a flat structure $\mathfrak{f}\in\mathcal{C}$ corresponds the unique 6-tuple  $(l_1,l_2,l_3,a_1,a_2,a_3)$ of $\mathcal{B}$ which determines $\mathfrak{f}$. Obviously we have $\Phi(\mathcal{C})=\mathcal{B}$.

\section{The Teichm\"uller space of flat pairs of pants with one singularity\label{teich}}

Let us denote by $\mathcal{F}(M)$ the space of all flat structures with one singularity on a pair of pants $M$. We fix an orientation on $M$ and let $\mathtt{Homeo}^+(M,\partial)$ be the set of homeomorphisms of $M$ which preserve the orientation and each boundary component of $M$ (setwise). It is well known that each element of $\mathtt{Homeo}^+(M,\partial)$ is isotopic to the identity (see Expos\'e 2 in \cite{0731.57001}). The space $\mathtt{Homeo}^+(M,\partial)$ acts on $\mathcal{F}(M)$ as follows: If $h\in \mathtt{Homeo}^+(M,\partial)$ and $\mathfrak{f}\in \mathcal{F}(M)$ then $(h,\mathfrak{f})\mapsto h*\mathfrak{f}$ where $h*\mathfrak{f}(x,y):=\mathfrak{f}(h(x),h(y))$.
\begin{defin}
\label{teichm}
 We define the Teichm\"uller space $\mathcal{T}(M)$ of $M$ as the quotient $\mathcal{F}(M)\slash \mathtt{Homeo}^+(M,\partial)$.
\end{defin}
Obviously, $\mathcal{T}(M)$ consists of all flat structures which belong to $\mathcal{C}$, the set of all flat structures defined by the parameters $l_i,\, a_i,\, i=1,2,3$. To each flat structure with one singularity on $M$ we may associate a unique configuration, which consists of a triangle $T$ and three rectangles $R_i,\,i=1,2,3$, glued as in Figure \ref{pantc1}. Conversely, to each configuration we may associate a unique flat structure with one singularity on $M$. This defines a mapping $\Phi:\mathcal{T}(M)\rightarrow \mathcal{B}$ which is one-to-one, and $\Phi(\mathcal{T}(M))=\mathcal{B}$.
\begin{rem}[Due to Athanase Papadopoulos]
 Let $\mathcal{T}_\mathcal{M}(M)$ be the Teichm\"uller space of flat structures with one singularity on a pair of pants defined by a pair of transverse measured foliations. This space is of dimension 4. Since $\mathcal{T}(M)$, the Teichm\"uller space of all flat structures with one singularity on a pair of pants, is of dimension 6, as we showed here, then $\mathcal{T}_\mathcal{M}(M)\subsetneq \mathcal{T}(M)$. Thus, the space of flat structures we counter here is larger than the space of flat structures induced by quadratic differentials with one zero. 
\end{rem}

\begin{prop}
 The mapping $\Phi:\mathcal{T}(M)\rightarrow \mathcal{B}$ is continuous and open. Thus, $\Phi$ is homeomorphism.
\end{prop}
\begin{proof}
 As explained in Proposition \ref{params2}, $\Phi$ is a bijection. If we consider the Euclidean distance between tuples of parameters $(l_1,l_2,l_3,a_1,a_2,a_3)$ and the following distance between two flat structures $\mathfrak{f}_1,\mathfrak{f}_2\in \mathcal{T}(M)$:
\[d(\mathfrak{f}_1,\mathfrak{f}_2)=\sup_{x,y\in M} |\mathfrak{f}_1(x,y)-\mathfrak{f}_2(x,y)|\]
we can see easily that both $\Phi$ and $\Phi^{-1}$ are continuous.
\end{proof}
Obviously, $\mathcal{B}$ is a non bounded convex subset of $\mathbb{R}^6$. In fact, the set of parameters $l_i,\,i=1,2,3$ can be seen as the space of all triangles which could be degenerate, but without those which have a side of length zero. This is, a non bounded pyramid of three faces with its sides and apex deleted. This space is convex in $\mathbb{R}^3$. The same can be said about the set of parameters $a_i,\,i=1,2,3$. Thus, $\mathcal{B}$ is convex in $\mathbb{R}^6$. 

The boundary $\partial\mathcal{B}$ has six  components which are non bounded convex subsets of $\mathbb{R}^5$. They correspond to cases of equality in the triangle inequality. Let us denote by $\partial_i \mathcal{B}_l\, ,\, i\in\{1,2,3\}$ the boundary component defined by
\[l_i = l_{(i+1)\mod 3} + l_{(i+2)\mod 3} \, .\]
The sets $\partial_i \mathcal{B}_l,\,i=1,2,3$ are pairwise disjoint subsets of $\partial \mathcal{B}$. Similarly, we denote by $\partial_j \mathcal{B}_a\, ,\, j\in\{1,2,3\}$ the boundary component defined by
\[a_j= a_{(j+1)\mod 3} + a_{(j+2)\mod 3} \, . \]
The sets $\partial_j \mathcal{B}_a,\,j=1,2,3$ are also pairwise disjoint subsets of $\partial \mathcal{B}$.

The intersection $\partial_{i,j}(\partial \mathcal{B)}=\partial_i\mathcal{B}_l \cap\partial_j\mathcal{B}_a$,  is homeomorphic to $\mathbb{R}^4$ if $i\neq j$, and empty if $i=j$. This means that each boundary component of $\mathcal{B}$ has two convex connected boundary components homeomorphic to $\mathbb{R}^4$. 

Then, by the homeomorphism $\Phi$, we have the following description of the Teichm\"uller space: 
\begin{cor}
 $\mathcal{T}(M)$ is homeomorphic to a non-compact submanifold of $\mathbb{R}^6$ of dimension $6$, with a natural cell-structure, having six cells of codimension one on its boundary.
\end{cor}
 We know that any convex subset of Euclidean n-space $\mathbb{E}^n$ is contractible. This gives us the following result:
\begin{thm}
 $\mathcal{T}(M)$ is a contractible space.
\end{thm}

In the next section we start a discussion of the decomposition into pairs of pants of closed surfaces with one singularity. Details will be given in subsequent work.
\section{Closed flat surfaces with one singularity \label{exmps}}
In general, a flat surface with one singularity is not decomposable by disjoint simple closed geodesics into pairs of pants. For this, we have the following example: 
\begin{exmp}
\label{nondecomp}
 Take a closed flat surface $M_3$ of genus $3$ with one singularity $s$. Suppose that we are able to decompose it by disjoint simple closed geodesics into pairs of pants. By the fact that there is no flat pair of pants without any singularity, we conclude that every pair of pants resulting from the decomposition has the singularity $s$ on its boundary. We know that a decomposition of $M_3$ by disjoint simple closed geodesics gives rise to four pairs of pants. This means that four boundary components should share the singularity. This is impossible under the present rules of identification (one boundary component is identified to one boundary component). Then, the decomposition is impossible.
\end{exmp}

The decomposition becomes possible if we change the rules of composition. For example, admitting geodesics to be non-simple or non-disjoint. That is, admitting the identification of parts of boundary components rather than boundary components entirely. I will not discuss here these rules, their study can be dealt with a separate work later.

The result in Example \ref{nondecomp} can be generalized as follows:

Let $M$ be a closed flat surface of genus $g\geq 3$ with a single conical
singularity $s$. Then $M$ cannot be decomposed into pairs of pants by disjoint simple closed geodesics.

\section*{Acknowledgements}
I would like to thank Athanase Papadopoulos for all the corrections and discussions during the preparation of this article. I also thank Daniele Alesandrini for his time. 

\bibliographystyle{abbrv}
\bibliography{bibliographarticle}
\textit{Author’s address:}\\
Ousama Malouf\\
Institut de Recherche Math\'ematique Avanc\'ee\\
7 rue Ren\'e Descartes\\
F-67084 Strasbourg Cedex\\
France\\
malouf@math.u-strasbg.fr
\end{document}